
%
%
%
\nonstopmode
\catcode`\@=11 
%
%
%
\font\seventeenrm=cmr17

\font\twelverm=cmr12
\font\ninerm=cmr9
\font\sixrm=cmr6

\font\seventeenbf=cmbx12 at 17pt
\font\fourteenbf=cmbx12 at 14pt
\font\twelvebf=cmbx12
\font\ninebf=cmbx9
\font\sixbf=cmbx6

\font\seventeeni=cmmi12 at 17pt             \skewchar\seventeeni='177
\font\fourteeni=cmmi12 at 14pt              \skewchar\fourteeni='177
\font\twelvei=cmmi12                        \skewchar\twelvei='177
\font\ninei=cmmi9                           \skewchar\ninei='177
\font\sixi=cmmi6                            \skewchar\sixi='177

\font\seventeensy=cmsy10 scaled\magstep3    \skewchar\seventeensy='60
\font\fourteensy=cmsy10 scaled\magstep2     \skewchar\fourteensy='60
\font\twelvesy=cmsy10 at 12pt               \skewchar\twelvesy='60
\font\ninesy=cmsy9                          \skewchar\ninesy='60
\font\sixsy=cmsy6                           \skewchar\sixsy='60

\font\seventeenex=cmex10 scaled\magstep3
\font\fourteenex=cmex10 scaled\magstep2
\font\twelveex=cmex10 at 12pt

\font\ninex=cmex10 at 9pt
\font\sevenex=cmex10 at 9pt
\font\sixex=cmex10 at 6pt
\font\fivex=cmex10 at 5pt

\font\seventeensl=cmsl10 scaled\magstep3
\font\fourteensl=cmsl10 scaled\magstep2
\font\twelvesl=cmsl10 scaled\magstep1
\font\ninesl=cmsl10 at 9pt
\font\sevensl=cmsl10 at 7pt
\font\sixsl=cmsl10 at 6pt
\font\fivesl=cmsl10 at 5pt

\font\seventeenit=cmti12 scaled\magstep2
\font\fourteenit=cmti12 scaled\magstep1
\font\twelveit=cmti12

\font\seventeentt=cmtt12 scaled\magstep2
\font\fourteentt=cmtt12 scaled\magstep1
\font\twelvett=cmtt12

\font\seventeencp=cmcsc10 scaled\magstep3
\font\fourteencp=cmcsc10 scaled\magstep2
\font\twelvecp=cmcsc10 scaled\magstep1
\font\tencp=cmcsc10

\newfam\cpfam

\font\seventeenss=cmss17
\font\fourteenss=cmss12 at 14pt
\font\twelvess=cmss12
\font\tenss=cmss10
\font\niness=cmss9

\font\sevenss=cmss8 at 7pt
\font\sixss=cmss8 at 6pt
\font\fivess=cmss8 at 5pt
\newfam\ssfam
\newdimen\b@gheight             \b@gheight=12pt
\newcount\f@ntkey               \f@ntkey=0
\def\f@m{\afterassignment\samef@nt\f@ntkey=}
\def\samef@nt{\fam=\f@ntkey \the\textfont\f@ntkey\relax}
\def\rm{\f@m0 }
\def\mit{\f@m1 }         
\def\cal{\f@m2 }
\def\it{\f@m\itfam}
\def\sl{\f@m\slfam}
\def\bf{\f@m\bffam}
\def\tt{\f@m\ttfam}
\def\ssf{\f@m\ssfam}
\def\caps{\f@m\cpfam}
\def\seventeenpoint{\relax
    \textfont0=\seventeenrm          \scriptfont0=\twelverm
      \scriptscriptfont0=\ninerm
    \textfont1=\seventeeni           \scriptfont1=\twelvei
      \scriptscriptfont1=\ninei
    \textfont2=\seventeensy          \scriptfont2=\twelvesy
      \scriptscriptfont2=\ninesy
    \textfont3=\seventeenex          \scriptfont3=\twelveex
      \scriptscriptfont3=\ninex
    \textfont\itfam=\seventeenit    
    \textfont\slfam=\seventeensl    
      \scriptscriptfont\slfam=\ninesl
    \textfont\bffam=\seventeenbf     \scriptfont\bffam=\twelvebf
      \scriptscriptfont\bffam=\ninebf
    \textfont\ttfam=\seventeentt
    \textfont\cpfam=\seventeencp
    \textfont\ssfam=\seventeenss     \scriptfont\ssfam=\twelvess
      \scriptscriptfont\ssfam=\niness
    \samef@nt
    \b@gheight=17pt
    \setbox\strutbox=\hbox{\vrule height 0.85\b@gheight
                                depth 0.35\b@gheight width\z@ }}
\def\fourteenpoint{\relax
    \textfont0=\fourteencp          \scriptfont0=\tenrm
      \scriptscriptfont0=\sevenrm
    \textfont1=\fourteeni           \scriptfont1=\teni
      \scriptscriptfont1=\seveni
    \textfont2=\fourteensy          \scriptfont2=\tensy
      \scriptscriptfont2=\sevensy
    \textfont3=\fourteenex          \scriptfont3=\twelveex
      \scriptscriptfont3=\tenex
    \textfont\itfam=\fourteenit     \scriptfont\itfam=\tenit
    \textfont\slfam=\fourteensl     \scriptfont\slfam=\tensl
      \scriptscriptfont\slfam=\sevensl
    \textfont\bffam=\fourteenbf     \scriptfont\bffam=\tenbf
      \scriptscriptfont\bffam=\sevenbf
    \textfont\ttfam=\fourteentt
    \textfont\cpfam=\fourteencp
    \textfont\ssfam=\fourteenss     \scriptfont\ssfam=\tenss
      \scriptscriptfont\ssfam=\sevenss
    \samef@nt
    \b@gheight=14pt
    \setbox\strutbox=\hbox{\vrule height 0.85\b@gheight
                                depth 0.35\b@gheight width\z@ }}
\def\twelvepoint{\relax
    \textfont0=\twelverm          \scriptfont0=\ninerm
      \scriptscriptfont0=\sixrm
    \textfont1=\twelvei           \scriptfont1=\ninei
      \scriptscriptfont1=\sixi
    \textfont2=\twelvesy           \scriptfont2=\ninesy
      \scriptscriptfont2=\sixsy
    \textfont3=\twelveex          \scriptfont3=\ninex
      \scriptscriptfont3=\sixex
    \textfont\itfam=\twelveit    
    \textfont\slfam=\twelvesl    
      \scriptscriptfont\slfam=\sixsl
    \textfont\bffam=\twelvebf     \scriptfont\bffam=\ninebf
      \scriptscriptfont\bffam=\sixbf
    \textfont\ttfam=\twelvett
    \textfont\cpfam=\twelvecp
    \textfont\ssfam=\twelvess     \scriptfont\ssfam=\niness
      \scriptscriptfont\ssfam=\sixss
    \samef@nt
    \b@gheight=12pt
    \setbox\strutbox=\hbox{\vrule height 0.85\b@gheight
                                depth 0.35\b@gheight width\z@ }}
\def\tenpoint{\relax
    \textfont0=\tenrm          \scriptfont0=\sevenrm
      \scriptscriptfont0=\fiverm
    \textfont1=\teni           \scriptfont1=\seveni
      \scriptscriptfont1=\fivei
    \textfont2=\tensy          \scriptfont2=\sevensy
      \scriptscriptfont2=\fivesy
    \textfont3=\tenex          \scriptfont3=\sevenex
      \scriptscriptfont3=\fivex
    \textfont\itfam=\tenit     \scriptfont\itfam=\seveni
    \textfont\slfam=\tensl     \scriptfont\slfam=\sevensl
      \scriptscriptfont\slfam=\fivesl
    \textfont\bffam=\tenbf     \scriptfont\bffam=\sevenbf
      \scriptscriptfont\bffam=\fivebf
    \textfont\ttfam=\tentt
    \textfont\cpfam=\tencp
    \textfont\ssfam=\tenss     \scriptfont\ssfam=\sevenss
      \scriptscriptfont\ssfam=\fivess
    \samef@nt
    \b@gheight=10pt
    \setbox\strutbox=\hbox{\vrule height 0.85\b@gheight
                                depth 0.35\b@gheight width\z@ }}
%
%
%
\normalbaselineskip = 15pt plus 0.2pt minus 0.1pt 
\normallineskip = 1.5pt plus 0.1pt minus 0.1pt
\normallineskiplimit = 1.5pt
\newskip\normaldisplayskip
\normaldisplayskip = 15pt plus 5pt minus 10pt 
\newskip\normaldispshortskip
\normaldispshortskip = 6pt plus 5pt
\newskip\normalparskip
\normalparskip = 6pt plus 2pt minus 1pt
\newskip\skipregister
\skipregister = 5pt plus 2pt minus 1.5pt
\newif\ifsingl@    \newif\ifdoubl@
\newif\iftwelv@    \twelv@true
\def\singlespace{\singl@true\doubl@false\spaces@t}
\def\doublespace{\singl@false\doubl@true\spaces@t}
\def\normalspace{\singl@false\doubl@false\spaces@t}
\def\Tenpoint{\tenpoint\twelv@false\spaces@t}
\def\Twelvepoint{\twelvepoint\twelv@true\spaces@t}
\def\spaces@t{\relax
      \iftwelv@ \ifsingl@\subspaces@t3:4;\else\subspaces@t1:1;\fi
       \else \ifsingl@\subspaces@t3:5;\else\subspaces@t4:5;\fi \fi
      \ifdoubl@ \multiply\baselineskip by 5
         \divide\baselineskip by 4 \fi }
\def\subspaces@t#1:#2;{
      \baselineskip = \normalbaselineskip
      \multiply\baselineskip by #1 \divide\baselineskip by #2
      \lineskip = \normallineskip
      \multiply\lineskip by #1 \divide\lineskip by #2
      \lineskiplimit = \normallineskiplimit
      \multiply\lineskiplimit by #1 \divide\lineskiplimit by #2
      \parskip = \normalparskip
      \multiply\parskip by #1 \divide\parskip by #2
      \abovedisplayskip = \normaldisplayskip
      \multiply\abovedisplayskip by #1 \divide\abovedisplayskip by #2
      \belowdisplayskip = \abovedisplayskip
      \abovedisplayshortskip = \normaldispshortskip
      \multiply\abovedisplayshortskip by #1
        \divide\abovedisplayshortskip by #2
      \belowdisplayshortskip = \abovedisplayshortskip
      \advance\belowdisplayshortskip by \belowdisplayskip
      \divide\belowdisplayshortskip by 2
      \smallskipamount = \skipregister
      \multiply\smallskipamount by #1 \divide\smallskipamount by #2
      \medskipamount = \smallskipamount \multiply\medskipamount by 2
      \bigskipamount = \smallskipamount \multiply\bigskipamount by 4 }
\def\normalbaselines{ \baselineskip=\normalbaselineskip
   \lineskip=\normallineskip \lineskiplimit=\normallineskip
   \iftwelv@\else \multiply\baselineskip by 4 \divide\baselineskip by 5
     \multiply\lineskiplimit by 4 \divide\lineskiplimit by 5
     \multiply\lineskip by 4 \divide\lineskip by 5 \fi }
\Twelvepoint  
%
\interlinepenalty=50
\interfootnotelinepenalty=5000
\predisplaypenalty=9000
\postdisplaypenalty=500
\hfuzz=1pt
\vfuzz=0.2pt
\dimen\footins=24 truecm 
\hoffset=10.5truemm 
\voffset=-8.5 truemm 
%
%
%
%
%
%
\def\footnote#1{\edef\@sf{\spacefactor\the\spacefactor}#1\@sf
      \insert\footins\bgroup\singl@true\doubl@false\Tenpoint
      \interlinepenalty=\interfootnotelinepenalty \let\par=\endgraf
        \leftskip=\z@skip \rightskip=\z@skip
        \splittopskip=10pt plus 1pt minus 1pt \floatingpenalty=20000
        \smallskip\item{#1}\bgroup\strut\aftergroup\@foot\let\next}
\skip\footins=\bigskipamount 
\dimen\footins=24truecm 
\newcount\fnotenumber
\def\clearfnotenumber{\fnotenumber=0}
\def\fnote{\advance\fnotenumber by1 \footnote{$^{\the\fnotenumber}$}}
\clearfnotenumber
%
%
\newcount\secnumber
\newcount\appnumber
\newif\ifs@c 
\newif\ifs@cd 
\s@cdtrue 
\def\unsectioned{\s@cdfalse\let\section=\subsection}
\def\clearappnumber{\appnumber=64}
\def\clearsecnumber{\secnumber=0}
\newskip\sectionskip         \sectionskip=\medskipamount
\newskip\headskip            \headskip=8pt plus 3pt minus 3pt
\newdimen\sectionminspace    \sectionminspace=10pc
\newdimen\referenceminspace  \referenceminspace=25pc
\def\Titlestyle#1{\par\begingroup \interlinepenalty=9999
     \leftskip=0.02\hsize plus 0.23\hsize minus 0.02\hsize
     \rightskip=\leftskip \parfillskip=0pt
     \advance\baselineskip by 0.5\baselineskip
     \hyphenpenalty=9000 \exhyphenpenalty=9000
     \tolerance=9999 \pretolerance=9000
     \spaceskip=0.333em \xspaceskip=0.5em
     \seventeenpoint
  \noindent #1\par\endgroup }
\def\titlestyle#1{\par\begingroup \interlinepenalty=9999
     \leftskip=0.02\hsize plus 0.23\hsize minus 0.02\hsize
     \rightskip=\leftskip \parfillskip=0pt
     \hyphenpenalty=9000 \exhyphenpenalty=9000
     \tolerance=9999 \pretolerance=9000
     \spaceskip=0.333em \xspaceskip=0.5em
     \fourteenpoint
   \noindent #1\par\endgroup }
%
\def\spacecheck#1{\dimen@=\pagegoal\advance\dimen@ by -\pagetotal
   \ifdim\dimen@<#1 \ifdim\dimen@>0pt \vfil\break \fi\fi}
\def\section#1{\cleareqnumber \s@ctrue \global\advance\secnumber by1
   \message{Section \the\secnumber: #1}
   \par \ifnum\the\lastpenalty=30000\else
   \penalty-200\vskip\sectionskip \spacecheck\sectionminspace\fi
   \noindent {\caps\enspace\S\the\secnumber\quad #1}\par
   \nobreak\vskip\headskip \penalty 30000 }
\def\subsection#1{\par
   \ifnum\the\lastpenalty=30000\else \penalty-100\smallskip
   \spacecheck\sectionminspace\fi
   \noindent\undertext{#1}\enspace \vadjust{\penalty5000}}

\def\undertext#1{\vtop{\hbox{#1}\kern 1pt \hrule}}
\def\subsubsection#1{\par
   \ifnum\the\lastpenalty=30000\else \penalty-100\smallskip \fi
   \noindent\hbox{#1}\enspace \vadjust{\penalty5000}}

\def\appendix#1{\cleareqnumber \s@cfalse \global\advance\appnumber by1
   \message{Appendix \char\the\appnumber: #1}
   \par \ifnum\the\lastpenalty=30000\else
   \penalty-200\vskip\sectionskip \spacecheck\sectionminspace\fi
   \noindent {\caps\enspace Appendix \char\the\appnumber\quad #1}\par
   \nobreak\vskip\headskip \penalty 30000 }
\clearsecnumber
\clearappnumber
%
%
\def\ack{\par\penalty-100\medskip \spacecheck\sectionminspace
   \line{\iftwelv@\fourteencp\else\twelvecp\fi\hfil {\caps
Acknowledgements}\hfil}%
\nobreak\vskip\headskip }
\def\refs{\begingroup \par\penalty-100\medskip \spacecheck\sectionminspace
   \line{\iftwelv@\fourteencp\else\twelvecp\fi\hfil REFERENCES\hfil}%
\nobreak\vskip\headskip \frenchspacing }
\def\endrefs{\par\endgroup}
%
\newcount\refnumber
\def\clearrefnumber{\refnumber=0}  \clearrefnumber
\newwrite\R@fs                              
\immediate\openout\R@fs=\jobname.references 
\def\closerefs{\immediate\closeout\R@fs} 
\def\refsout{\closerefs\refs
\catcode`\@=11                          
\input\jobname.references               
\catcode`\@=12			        
\endrefs}
\def\refitem#1{\item{{\bf #1}}}
\def\ifundefined#1{\expandafter\ifx\csname#1\endcsname\relax}
%
%
\def\[#1]{\ifundefined{#1R@FNO}%
\global\advance\refnumber by1%
\expandafter\xdef\csname#1R@FNO\endcsname{[\the\refnumber]}%
\immediate\write\R@fs{\noexpand\refitem{\csname#1R@FNO\endcsname}%
\noexpand\csname#1R@F\endcsname}\fi{\bf \csname#1R@FNO\endcsname}}
\def\refdef[#1]#2{\expandafter\gdef\csname#1R@F\endcsname{{#2}}}
%
%
%
%
%
%
\newcount\eqnumber
\def\cleareqnumber{\eqnumber=0}
\newif\ifal@gn \al@gnfalse  
\def\veqnalign#1{\al@gntrue \vbox{\eqalignno{#1}} \al@gnfalse}
\def\eqnalign#1{\al@gntrue \eqalignno{#1} \al@gnfalse}
\def\(#1){\relax%
\ifundefined{#1@Q}
 \global\advance\eqnumber by1
 \ifs@cd
  \ifs@c
   \expandafter\xdef\csname#1@Q\endcsname{{%
\noexpand\rm(\the\secnumber .\the\eqnumber)}}
  \else
   \expandafter\xdef\csname#1@Q\endcsname{{%
\noexpand\rm(\char\the\appnumber .\the\eqnumber)}}
  \fi
 \else
  \expandafter\xdef\csname#1@Q\endcsname{{\noexpand\rm(\the\eqnumber)}}
 \fi
 \ifal@gn
    & \csname#1@Q\endcsname
 \else
    \eqno \csname#1@Q\endcsname
 \fi
\else%
\csname#1@Q\endcsname\fi\global\let\@Q=\relax}
%
%
\newif\iffrontpage \frontpagefalse
\headline={\hfil}
\footline={\iffrontpage\hfil\else \hss\twelverm
-- \folio\ --\hss \fi }
\def\monthname{\relax\ifcase\month 0/\or January\or February\or
   March\or April\or May\or June\or July\or August\or September\or
   October\or November\or December\else\number\month/\fi}
\hsize=14 truecm
\vsize=22 truecm
\skip\footins=\bigskipamount
\normalspace
%
%
%
\newskip\frontpageskip
\newif\ifp@bblock \p@bblocktrue
\newif\ifm@nth \m@nthtrue
\newtoks\pubnum
\newtoks\pubtype
\newtoks\m@nthn@me
\newcount\Ye@r
\advance\Ye@r by \year
\advance\Ye@r by -1900
\def\Year#1{\Ye@r=#1}
\def\Month#1{\m@nthfalse \m@nthn@me={#1}}
\def\m@nthname{\ifm@nth\monthname\else\the\m@nthn@me\fi}
\def\titlepage{\global\frontpagetrue\hrule height\z@ \relax
               \ifp@bblock\pubblock\fi\relax }
\def\endtitlepage{\vfil\break
                  \frontpagefalse} 
\def\bonntitlepage{\global\frontpagetrue\hrule height\z@ \relax
               \ifp@bblock\pubblock\fi\relax }
\frontpageskip=12pt plus .5fil minus 2pt
\pubtype={\iftwelv@\twelvesl\else\tensl\fi\ (Preliminary Version)}
\pubnum={?}
\def\nopubblock{\p@bblockfalse}
\def\pubblock{\line{\hfil\iftwelv@\twelverm\else\tenrm\fi%
US--\number\Ye@r--\the\pubnum\the\pubtype}
              \line{\hfil\iftwelv@\twelverm\else\tenrm\fi%
\m@nthname\ \number\year}}
\def\title#1{\vskip\frontpageskip\Titlestyle{\caps #1}\vskip3\headskip}
\def\author#1{\vskip.5\frontpageskip\titlestyle{\caps #1}\nobreak}
\def\andauthor{\vskip.5\frontpageskip\centerline{and}\author}

\def\address#1{\par\kern 5pt\titlestyle{
\it #1}}
\def\andaddress{\par\kern 5pt \centerline{\sl and} \address}

\def\Santiago{\address{Departamento de F{\'\i}sica de Part{\'\i}culas
Elementales\break
Universidad de Santiago, Santiago de Compostela 15706, SPAIN}}
\def\Montpellier{\address{Laboratoire de Physique Math\'ematique\break
Universit\'e de Montpellier II,Place Eug\`ene Bataillon\break
34095 Montpellier, CEDEX 5, FRANCE}}

\def\abstract#1{\par\dimen@=\prevdepth \hrule height\z@ \prevdepth=\dimen@
   \vskip\frontpageskip\spacecheck\sectionminspace
   \centerline{\iftwelv@\fourteencp\else\twelvecp\fi ABSTRACT}\vskip\headskip
   {\noindent #1}}
%

%
%
%
\def\leaderfill{\leaders\hbox to 1em{\hss.\hss}\hfill}
\def\boxit#1{\vcenter{\hrule\hbox{\vrule\kern8pt
      \vbox{\kern8pt#1\kern8pt}\kern8pt\vrule}\hrule}}

%
%
%
\def\ref#1{{\bf [#1]}}
\def\ie{{\it i.e.\/}}
\def\nl{\hfil\break}
%
%
%
%
%
\newif\ifm@thstyle \m@thstylefalse
\def\mathstyle{\m@thstyletrue}
\def\proclaim#1#2\par{\smallbreak\begingroup
\advance\baselineskip by -0.25\baselineskip%
\advance\belowdisplayskip by -0.35\belowdisplayskip%
\advance\abovedisplayskip by -0.35\abovedisplayskip%
    \noindent{\caps#1.\enspace}{#2}\par\endgroup%
\smallbreak}
\def\m@kem@th<#1>#2#3{%
\ifm@thstyle \global\advance\eqnumber by1
 \ifs@cd
  \ifs@c
   \expandafter\xdef\csname#1\endcsname{{%
\noexpand #2\ \the\secnumber .\the\eqnumber}}
  \else
   \expandafter\xdef\csname#1\endcsname{{%
\noexpand #2\ \char\the\appnumber .\the\eqnumber}}
  \fi
 \else
  \expandafter\xdef\csname#1\endcsname{{\noexpand #2\ \the\eqnumber}}
 \fi
 \proclaim{\csname#1\endcsname}{#3}
\else
 \proclaim{#2}{#3}
\fi}
%
%
%
%
%
%
\def\Thm<#1>#2{\m@kem@th<#1M@TH>{Theorem}{\sl#2}}
\def\Prop<#1>#2{\m@kem@th<#1M@TH>{Proposition}{\sl#2}}
\def\Def<#1>#2{\m@kem@th<#1M@TH>{Definition}{\rm#2}}
\def\Lem<#1>#2{\m@kem@th<#1M@TH>{Lemma}{\sl#2}}
\def\Cor<#1>#2{\m@kem@th<#1M@TH>{Corollary}{\sl#2}}
\def\Conj<#1>#2{\m@kem@th<#1M@TH>{Conjecture}{\sl#2}}
\def\Rmk<#1>#2{\m@kem@th<#1M@TH>{Remark}{\rm#2}}
\def\Exm<#1>#2{\m@kem@th<#1M@TH>{Example}{\rm#2}}
\def\Qry<#1>#2{\m@kem@th<#1M@TH>{Query}{\it#2}}
\def\Example<#1>#2{\m@kem@th<#1M@TH>{Example}{\it#2}}
\def\Exercise<#1>#2{\m@kem@th<#1M@TH>{Execise}{\it#2}}
\def\Proof{\noindent{\caps Proof:}\enspace}

\let\Example=\Exm
\def\<#1>{\csname#1M@TH\endcsname}
%
%
\def\qed{\vrule width 0.6em height 0.5em depth 0.2em}

\def\lapprox{\hbox{\lower3pt\hbox{$\buildrel<\over\sim$}}}
\def\gapprox{\hbox{\lower3pt\hbox{$\buildrel<\over\sim$}}}
\def\quotient#1#2{#1/\lower0pt\hbox{${#2}$}}
%
%
\def\to{\rightarrow}
%

%
%
%
\def\nats{{\bf N}} 
%
%
\def\Tr{{\rm Tr}}
\def\underrightarrow#1{\vtop{\ialign{##\crcr
      $\hfil\displaystyle{#1}\hfil$\crcr
      \noalign{\kern-\p@\nointerlineskip}
      \rightarrowfill\crcr}}} 
\def\underleftarrow#1{\vtop{\ialign{##\crcr
      $\hfil\displaystyle{#1}\hfil$\crcr
      \noalign{\kern-\p@\nointerlineskip}
      \leftarrowfill\crcr}}}  

%
%
\def\comm#1#2{\left[#1\, ,\,#2\right]}
\def\anticomm#1#2{\left\{#1\, ,\,#2\right\}}
%
%
\def\pder#1#2{{{\partial #1}\over{\partial #2}}}
%
%
%
%
%
%
\newdimen\unit
\newdimen\redunit
%
%
\def\p@int#1:#2 #3 {\rlap{\kern#2\unit
     \raise#3\unit\hbox{#1}}}
%
%
\def\th@r{\vrule height0\unit depth.1\unit width1\unit}
\def\bh@r{\vrule height.1\unit depth0\unit width1\unit}
\def\lv@r{\vrule height1\unit depth0\unit width.1\unit}
\def\rv@r{\vrule height1\unit depth0\unit width.1\unit}
%
%
\def\t@ble@u{\hbox{\p@int\bh@r:0 0
                   \p@int\lv@r:0 0
                   \p@int\rv@r:.9 0
                   \p@int\th@r:0 1
                   }
             }
%
%
\def\t@bleau#1#2{\rlap{\kern#1\redunit
     \raise#2\redunit\t@ble@u}}
%
%
\newcount\n
\newcount\m
\def\makecol#1#2#3{\n=0 \m=#3
  \loop\ifnum\n<#1{}\advance\m by -1 \t@bleau{#2}{\number\m}\advance\n by 1
\repeat}
%
%
\def\makerow#1#2#3{\n=0 \m=#3
 \loop\ifnum\n<#1{}\advance\m by 1 \t@bleau{\number\m}{#2}\advance\n by 1
\repeat}
%
%
\def\checkunits{\ifinner \unit=6pt \else \unit=8pt \fi
                \redunit=0.9\unit } 
\def\ytsym#1{\checkunits\kern-.5\unit
  \vcenter{\hbox{\makerow{#1}{0}{0}\kern#1\unit}}\kern.5em} 
\def\ytant#1{\checkunits\kern.5em
  \vcenter{\hbox{\makecol{#1}{0}{0}\kern1\unit}}\kern.5em} 
\def\ytwo#1#2{\checkunits
  \vcenter{\hbox{\makecol{#1}{0}{0}\makecol{#2}{1}{0}\kern2\unit}}
                  \ } 
\def\ythree#1#2#3{\checkunits
  \vcenter{\hbox{\makecol{#1}{0}{0}\makecol{#2}{1}{0}\makecol{#3}{2}{0}%
\kern3\unit}}
                  \ } 

%
%
%
\def\PRL#1#2#3{{\sl Phys. Rev. Lett.} {\bf#1} (#2) #3}
\def\NPB#1#2#3{{\sl Nucl. Phys.} {\bf B#1} (#2) #3}

\def\CMP#1#2#3{{\sl Comm. Math. Phys.} {\bf #1} (#2) #3}

\def\PLB#1#2#3{{\sl Phys. Lett.} {\bf #1B} (#2) #3}
\def\JMP#1#2#3{{\sl J. Math. Phys.} {\bf #1} (#2) #3}

\def\Invm#1#2#3{{\sl Invent. math.} {\bf #1} (#2) #3}
\def\LMP#1#2#3{{\sl Letters in Math. Phys.} {\bf #1} (#2) #3}
\def\IJMPA#1#2#3{{\sl Int. J. Mod. Phys.} {\bf A#1} (#2) #3}

\def\TMP#1#2#3{{\sl Theor. Mat. Phys.} {\bf #1} (#2) #3}

\def\MPLA#1#2#3{{\sl Mod. Phys. Lett.} {\bf A#1} (#2) #3}

\def\PJAS#1#2#3{{\sl Proc. Jpn. Acad. Sci.} {\bf #1} (#2) #3}
\def\JPSJ#1#2#3{{\sl J. Phys. Soc. Jpn.} {\bf #1} (#2) #3}

\catcode`\@=12 
%
%
%


\def\d{\partial}
\def\pdo{{\hbox{$\Psi$DO}}}

\def\comb[#1/#2]{\left[{#1\atop#2}\right]}

\def\rflecha#1{\setbox1=\hbox{\ninerm ~#1~}%
\buildrel\hbox{\ninerm #1}\over{\hbox to\wd1{\rightarrowfill}}}
\def\lflecha#1{\setbox1=\hbox{\ninerm ~#1~}%
\buildrel\hbox{\ninerm #1}\over{\hbox to\wd1{\leftarrowfill}}}

\def\gradL#1{{\delta #1\ov \delta L}}
\def\gradtL#1{{\delta #1\ov \delta \tilde L}}
\def\gradA#1{{\delta #1\ov \delta A}}
\def\gradtA#1{{\delta #1\ov \delta \tilde A}}
\def\gradB#1{{\delta #1\ov \delta B}}

\def\ov{\over}
\def\fr#1/#2{\mathord{\hbox{${#1}\over{#2}$}}}

\def\GD{{\ssf GD}}
\def\W{\mathord{\ssf W}}

\def\ope[#1][#2]{{{#2}\over{\ifnum#1=1 {z-w} \else {(z-w)^{#1}}\fi}}}



\refdef[WReview]{P.~Bouwknegt and K.~Schoutens, {\sl ${\cal
W}$-Symmetry in Conformal Field Theory},  {\it Phys. Reps.} to
appear.}
\refdef[Univ]{J.~M.~Figueroa-O'Farrill and E.~Ramos,
\JMP{33}{1992}{833}.}
\refdef[Wn]{A.~B.~Zamolodchikov, \TMP{65}{1986}{1205};\nl
V.~A.~Fateev and S.~L.~Lykyanov, \IJMPA{3}{1988}{507}.}
\refdef[Dickey]{L.~A.~Dickey,  {\sl Soliton equations and Hamiltonian
systems},  Advanced Series in Mathematical Physics Vol.12,  World
Scientific Publ.~Co..}
\refdef[GD]{I.~M.~Gel'fand and L.~A.~Dickey, {\sl A family of
Hamiltonian structures connected with integrable nonlinear
differential equations}, Preprint 136, IPM AN SSSR, Moscow (1978).}
\refdef[Adler]{M.~Adler, \Invm{50}{1979}{403}.}
\refdef[KP]{E.~Date, M.~Jimbo, M.~Kashiwara, and T.~Miwa
\PJAS{57A}{1981}{387}; \JPSJ{50}{1981}{3866}.}
\refdef[WKP]{L.~A.~Dickey, {\sl Annals NY Acad.~Sci.} {\bf 491}(1987)
131;\nl J.~M.~Figueroa-O'Farrill, J.~Mas, and E.~Ramos,
\PLB{266}{1991}{298};\nl F.~Yu and Y.-S.~Wu, \NPB{373}{1992}{713}.}
\refdef[WKPq]{J.~M.~Figueroa-O'Farrill, J.~Mas, and E.~Ramos
\CMP{158}{1993}{17}.}
\refdef[WinftyKP]{K.~Yamagishi, \PLB{259}{1991}{436};\nl F.~Yu and
Y.-S.~Wu, \PLB{236}{1991}{220}}
\refdef[Class]{J.~M.~Figueroa-O'Farrill and E.~Ramos, \PLB{282}{1992}{357}
({\tt hep-th/9202040}).}
\refdef[Radul]{A.~O.~Radul, in {\sl Applied methods of nonlinear
analysis and control}, pp. 149-157, Mironov, Moroz, and Tshernjatin,
eds.,  MGU 1987 (Russian).}
\refdef[winfty]{I.~Bakas, \PLB{228}{1989}{406}; \CMP{134}{1990}{487}.}
\refdef[Winfty]{C.~N.~Pope, L.~J.~Romans, and X.~Shen,
\NPB{339}{1990}{191}.}
\refdef[Woneplusinfty]{C. N. Pope, L. J. Romans, and X. Shen,
\PLB{242}{1990}{401}.}
\refdef[matrix]{ L.~Bonora and C.S.~Xiong, \NPB{405}{1993}{191}.}
\refdef[redukp]{W. Oevel and W. Strampp,~  \CMP{157}{1993}{51-81}.\nl
H.~Aratyn, L.~A.~Ferreira, J.~F.~Gomes, and
A.~H.~Zimerman, ~\NPB{402}{1993}{85-117} \nl
L.~Bonora, Q.P.~Liu and C.S.~ Xiong,~ {\sl The
integrable hierarchy constructed from a pair of KdV type
hierarchies and its
associated $\W$ algebra}, ({\tt hep-th/9408035}); \nl
D. Depireux and J. Shiff,{\sl On the Hamiltonian structures
and the reductions of KP hierarchy},{\tt hep-th/9210080}. \nl
F. Toppan, {\sl On Coset Reductions of KP Hierarchies},
 {\tt hep-th/9409126. } }
\refdef[BaYu]{ I.~ Bakas and E.~ Kiritsis, \IJMPA{7}{1992}{55}, \nl
F.~Yu and Y.-S.~Wu, \PRL{68}{1992}{2996} ({\tt
hep-th/9112009}); \NPB{373}{1992}{713} }
\refdef[Yu]{ F. Yu, \LMP{29}{1993}{175}.}
\refdef[dickeyred]{L.~A.~Dickey,
{\sl On the constrained KP hierarchy, I and II}, ( {\tt hep-th/9407038} and
{\tt hep-th/9411005}).}
\refdef[oevel]{W. Oevel and W. Strampp,~  \CMP{157}{1993}{51-81}.}
\refdef[hidi]{F. Mart\'\i nez-Mor\'as and E. Ramos,~
 \CMP{157}{1993}{573-589}.
\nl
F. Mart\'\i nez-Mor\'as, J. Mas and E. Ramos,
{}~\MPLA{8}{1993}{2189-2197}.}
\refdef[kupwilson]{B.A. ~Kupershmidt and G.~ Wilson,~ \Invm{62}{1981}{403}.}
\refdef[Fateev]{V.A. Fateev and S.L. Lykyanov, \IJMPA{3}{1988}{507}.}
\def\pubblock{ \line{\hfil\tt USC-FT-1/95}
               \line{\hfil\tt PM/95-03}
               \line{\hfil\tt q-alg/9501009}
               \line{\hfil\rm January 1995}}
\mathstyle
\titlepage
\title{The Constrained KP Hierarchy and the Generalised Miura Transformation}
\author{Javier Mas \footnote{$^\sharp$}{\tt e-mail:
jamas@gaes.usc.es}}
\Santiago
\andauthor{Eduardo Ramos\footnote{$^\flat$}{\tt e-mail:
ramos@lpm.univ-montp2.fr}}
\Montpellier

\unsectioned
\abstract{ Recently much attention has been paid
to the restriction of KP to
the submanifold of operators which can be represented
as a ratio of two purely
differential operators $L=AB^{-1}$. Whereas most
of the aspects concerning this reduced
hierarchy,
like the Lax flows and the Hamiltonians, are by now
well understood, there still
lacks a clear and conclusive statement about the associated
Poisson structure.
We fill this gap by
placing the problem
in a more general framework and then showing how the
required result follows from
an interesting property of the second Gelfand-Dickey brackets under
multiplication and inversion of Lax operators. As a byproduct
we give an elegant and simple proof of the generalised
Kupershmidt-Wilson theorem.}
\endtitlepage

$\W$ algebras provide a link between two dimensional conformal field
theories and integrable systems. At the heart of their relation with the
former lies the so called ``free field realization", which in the context
of the latter receives the name of (quantized) Miura transformation
\[Fateev].
On the other hand, the approach based on integrable systems has provided
a unifying framework for those algebras via the formalism of
pseudo-differential operators (\pdo); {\ie} there exists a basic common
ingredient, the so called  Adler map, from which the $\W$ algebras can be
readily constructed as Poisson bracket algebras. The data that specify
a particular $\W$ algebra is encoded in the particular form
of the associated Lax operator.
This scheme has raised the hope of establishing a full classification
or atlas. Whereas only partial results have been achieved, the
increasing size of this atlas, with new examples being constructed everyday,
calls for a deeper understanding of the way in which different $\W$ algebras
can be related.

This letter is a modest step towards this direction.
Although many of the results presented here were already known, our
approach makes emphasis on the
astonishing simplicity that lies behind a number of important and, up to
now, disperse results, and whose individual proof required  a considerable
amount of insight and calculational thrust.

Our original motivation stemmed from the construction of free (multi)
boson realization for non-linear $\W_\infty$ type of
algebras\[BaYu].
 From the point of view of integrable systems, the associated piece of data
is a KP-Lax operator of the form $L = AB^{-1}$, where $A=(\d+\varphi_1)\cdots
(\d+\varphi_m)$ and $B=(\d+\phi_1)\cdots(\d+\phi_n)$. This kind of Lax operator
arises naturally in the context of matrix models \[matrix] and the
associated hierarchy and Hamiltonian structure has been the subject of
recent intense research \[redukp].
The possibility of inducing a $\W_\infty$ type of algebra from the free
boson Poisson brackets for $\varphi_i$ and $\phi_j$ is the content of what we
call the generalized Kupershmidt-Wilson theorem. As the reader will see,
this theorem follows as a trivial corollary from the property of
``self-similarity" of the second Gelfand-Dickey brackets against product and
inversion of Lax operators.

The natural geometrical arena for KP is the
space of pseudodifferential operators (\pdo ) of the form
$$
L =\d^n + \sum_{i\geq 1}u_i \d^{n-i}. \(laxnkp)
$$

Each $L$ parameterizes a point on
the infinite dimensional manifold ${\cal M}_L$ whose coordinates are the
functions in the
unbounded set $\{u_1,u_2,...\}$, which in addition may be considered as
generators
of a differential polynomial algebra ${\ssf A}_L$.
On this manifold we may consider functions or rather functionals $F(L)$,
$G(L)$ and Poisson
brackets thereof
$\anticomm{F}{G}$. The interesting Poisson brackets,
from the point of view of CFT, are the ones of the associated
integrable hierarchy, {\ie} the KP-hierarchy. They are called
{\sl second} Gelfand-Dickey brackets \[GD]\[Dickey], and are explicitly
given by
$$
\anticomm{F}{G} = \Tr \left(\gradL F(L\gradL G)_+L
-\gradL F L(\gradL G L)_+ \right).
\(geldidos)
$$
As usual, $\Tr~ L$ stands for $\int res ~L$, where
$res~ L$ is the coefficient of $\d^{-1}$ in $L$,
the subscript $+$ stands for the projection onto the differential
part, and
$\gradL F$ and $\gradL G$ are the gradients of $F$ and $G$
evaluated at the point $L$:
$$
\gradL F = \sum_{i\in \nats} \d^{-n+i-1}{\delta F\ov \delta u_i}.
\(gradients)
$$
It is a nontrivial fact that these formulas actually define
consistent Poisson brackets. In fact,
an essential point behind this property is
the form of the Lax operator $L$.  If $L$ did not
have an infinite tail or, in
general, if some constraints appeared among the fields $u_i$ the bracket
as defined above
would not enjoy the Poisson properties, unless the definition of the
gradients would be modified accordingly, an extremely difficult task
in general. However there exist some
peculiar constrained forms for $L$ for which the formula  \(geldidos)
still defines consistent Poisson brackets as it stands
without further modification. The first well-known case is when $L$ is a
purely differential operator, \ie $~u_{n+1}=u_{n+2}=....=0$.
This reduction leads to the Poisson structures associated with the famous
 n-KdV
hierarchy.

We will call Lax operators for which \(geldidos) defines a consistent
set of Poisson
brackets, {\sl consistent Lax operators}, of which we have already given
two examples: the so called
unreduced n-KP, where $n\in {\bf Z},~ (n\neq 0)$,
and n-KdV Lax
 operators with
$n\in {\bf N}$.
The expert in constrained dynamical systems will readily recognise
that, in this context, {\sl consistent} is the same as either unconstrained,
or else, first class constrained.

{}From \(geldidos),
a set of fundamental Poisson brackets among the generators $\{u_i\}$
of the polinomial
algebra ${\ssf A}_L$ may be derived in a straightforward manner by considering
linear functionals.
For $L$ the n-KP operator in \(laxnkp), the infinite
set $u_1,u_2,...$  spans a non-linear algebra  called
$\W_{KP}^{(n)}$ \[Radul]\[WKPq] (also named ${\hat{\W}}_{1+\infty}^{(n)}$).
When $L$ is the n-KdV operator, after setting  $u_{n+1}=u_{n+2}=....=0$ in
$\W_{KP}^{(n)}$ the remaining
fields close another
non-linear
algebra named $\GD_n$ after I.M. Gelfand and L.A. Dickey.

All the information about the Poisson brackets is neatly
encoded in the essential ``geometrical" structure
behind it: the Adler map. For any {\sl one-form }
$X= \sum_i \d^{-i+n-1}x_i$, dual
to the {\sl deformations} of $L$ under the pairing defined by the
trace, the Adler map is then defined as
$$
J_L(X) = (LX)_+L - L(XL)_+ \(adler)
$$
Therefore, the
Poisson brackets in \(geldidos) can be equivalently written as
$$
\anticomm{F}{G} = \Tr ~\gradL F J_L\left(\gradL G \right).\(compois)
$$

The question
 of the consistency of $L$ can be traced back to the fact that
$J_L(X)$ should be understood as a deformation
of $L$. Hence the Adler map, in order not take us
out of the (constrained) manifold, should yield
an operator of the same form as $L-\d^n$.

{\caps Theorem I}:
{\sl Let $A = \sum_i a_i\d^i$ and $B = \sum_j b_j \d^j$ be
two consistent
Lax \pdo, and let the Poisson brackets among the generators
$\{a_i\}$ of ${\ssf A}_A$ and $\{b_j\}$ of ${\ssf A}_B$ be given
by the second Gelfand-Dickey bracket
\(geldidos); moreover let  $\anticomm{a_i}{b_j} = 0~~
\forall ~i,j$. Then we can form the product $AB = L = \sum_k
u_k\d^k $ and the Poisson brackets induced on the  $u_k$'s
via the embedding
$u_k = u_k(a_i,b_j)$ are nothing but the second Gelfand-Dickey
brackets computed
on ${\ssf A}_L$; L is again a consistent Lax operator.}

\Proof
Considered as functions of $A$ and $B$ the functional $F$ has
infinitesimal increment of the form
$$
\delta F(A,B) = \Tr~ (\gradA F\delta A +\gradB F\delta B) \(equno)
$$
whereas considered as a function of $L$ the same increment
looks like
$$
\delta F(L) = \Tr~ \gradL F \delta L =
\Tr ~\gradL F (\delta A~ B +A~\delta B)
$$
using the cyclic property of the trace and comparing with \(equno)
we arrive at the following identification
$$
\gradA F = B\gradL F~~;~~
\gradB F = \gradL F A. \(gradientes)
$$
The computation is now straightforward
$$
\veqnalign{
\anticomm{F}{G} & = \cr
=&~ \Tr ~\left(
  \gradA F(A\gradA G)_+A -\gradA F A(\gradA G A)_+  ~+~
  \gradB F(B\gradB G)_+B -\gradB F B(\gradB G B)_+  \right) \cr
=&~ \Tr ~B\gradL F(AB \gradL G)_+A -
 \Tr~ B\gradL F A(B\gradL F A)_+ \cr
& + \Tr ~\gradL F A(B\gradL G A)_+B -
\Tr~\gradL F AB(\gradL G AB)_+ \cr
=&~ \Tr~\left( L\gradL F(L \gradL G)_+ -
\gradL F L(\gradL G L)_+ \right).
\(prueba1) }
$$
\hfill\qed

The carefull reader may suspect that we have sloppily omitted
projectors onto differential or integral parts.
Namely, suppose that A were
purely differential, then the right hand of both
equalities in \(gradientes)
should be enclosed in $(~)_-$. However this is not necessary because,
as usual, the Adler map in the expression for
the Poisson brackets takes care of all projections automatically.

We could have phrased Theorem I in a different way
which emphasizes
the fundamental Poisson-bracket algebra. For example, if $A$
and $B$ are KdV (differential) operators of order $n$ and $m$,
the associated
non-linear algebras spanned by the $\{a_i\}$ and the $\{b_j\}$
are $\GD_n$ and $\GD_m$. As a result
of this theorem, we find that the set of $n+m$ differential
polynomials $\{u_k =
u_k(a_i,b_j)\}$ defined by the relation $L=AB$ span
the algebra $\GD_{n+m}$. Hence
$\GD_{n+m}$ is a Poisson subalgebra of $\GD_n\times\GD_m$.

{\caps Corollary I}: {\sl  The Kupershmidt-Wilson
theorem \[kupwilson] follows;
namely, the
Miura transformation defines a Poisson algebra homomorphism.}

\Proof
First, notice that the content of theorem I admits,
by repeated application,
a straightforward extension to
the case of multiple factorization $L = ABC\cdots D$.
Next choose $A = \d+a$ and $B = \d+b$. The specialization of the
first expression in \(prueba1) to these cases leads to
$$
\veqnalign{
\anticomm{F}{G}
& = \int \left( {\delta F\ov \delta a}\right)
\left({\delta G\ov \delta a}\right)'
+\left( {\delta F\ov \delta b}\right)\left({\delta G\ov
\delta b}\right)'
\cr}
$$
\ie, the basic building blocks  $a$ and $b$ are mutually
commuting free fields:
$$
\anticomm{a(x)}{a(y)} = \delta'(x-y)~~~~~
\anticomm{b(x)}{b(y)} = \delta'(x-y)~~~~~
\anticomm{a(x)}{b(y)} = 0
$$
\hfill\qed

 {\caps Theorem II}: {\sl Let $A=\sum_i a_i\d^i$ be an invertible
consistent   Lax operator and let the Poisson brackets among
the $a_i$ be given
by \(geldidos). Take  $L=A^{-1} = \sum_j u_j\d^j$. The
Poisson brackets among the $u_j$ induced by the
mapping $u_j \to u_j(a_i)$  are
given again by the second Gelfand-Dickey bracket on ${\ssf A}_L$
except for a
relative minus sign; $L$ is  also a consistent Lax operator.}

\Proof
Following {\sl mutatis mutandi}, identical steps as for Theorem I,
we have on one side
$$
\gradA F =- L \gradL F L~,
$$
and, by direct substitution we find
$$
\veqnalign{
\anticomm{F}{G}(A) &=
\Tr ~\left(
  \gradA F(A\gradA G)_+A -\gradA F A(\gradA G A)_+\right) \cr
&=
-\Tr ~\left(
  \gradL F(L\gradL G)_+L -\gradL FL(\gradL G L)_+ \right)  \cr
&= - \anticomm{F}{G}(L) }
$$
\hfill \qed

We consider both theorems as a kind of generating rules
for consistent Lax operators
under the Adler map. Of special importance for KP is the second
one since it constitutes a
kind of Miura transform for the KP operator. Namely writing
$L = A^{-1}$ where
$A= \d+a$ the potentials in $L = \d^{-1} + u_1\d^{-2} + u_2\d^{-3} +...$
are expressed as differential polynomials in $a$.
The algebra of the former ones is
named $\W^{(-1)}_{KP}$ \[Radul]\[WKPq], and is induced by
the embedding $u_i\rightarrow u_i(a)$
where $a$ satisfies the free field Poisson brackets. More
generally, if $A$ is of order $n$ and its fields $a_i$
span the algebra $\GD_n$,
$L$ is of order $-n$ and its fields span the algebra $\W_{KP}^{(-n)}$

The recent literature on the topic is concerned with the so called {\sl
multi-boson} reduction of KP. In particular the
juxtaposition of theorems
I and II helps to clarify the analysis of
ref. \[dickeyred] in what concerns the Poisson structure
induced by the mapping
$A,B\to AB^{-1}$. In particular it explains the fact that
the Lax flows are
hamiltonian with respect to the direct ``difference" of hamiltonian
structures for $A$ and $B$.

The two comments above can be more precisely stated in the form of two
corollaries as follows:

{\caps Corollary II}:
{\sl Let $A = \sum_i a_i\d^i$ and $B = \sum_j b_j \d^j$ be
two consistent
Lax \pdo, and let the Poisson brackets among the generators
$\{a_i\}$ of ${\ssf A}_A$ and $\{b_j\}$ of ${\ssf A}_B$ be given
respectively
by the Gelfand-Dickey bracket \(geldidos) with opposite signs; moreover let
$\anticomm{a_i}{b_j} = 0~~ \forall ~i,j$. Then we can form the product
$AB^{-1} = L = \sum_k u_k\d^k $ and the Poisson brackets induced on the
$u_k$'s via the embedding
$u_k = u_k(a_i,b_j)$ are equal to the second Gelfand-Dickey
brackets computed
on ${\ssf A}_L$; L is again a consistent Lax operator.}

{\caps Corollary III}: {\sl Kupershmidt-Wilson-Yu
theorem\[kupwilson]\[Yu]. Let $L$
be a {\pdo}  of order $n-m$.
If $L$ admits a factorization (generalised Miura transformation)
$L= (\d + \varphi_1)\cdots (\d +\varphi_n)(\d + \phi_1)^{-1}
\cdots (\d +\phi_m)^{-1}=\sum_k u_k\d^k$, and the only nonzero
Poisson brackets among the generalised Miura fields are given
by $\{\varphi_i(y),\varphi_j(x)\}=\delta_{ij}\delta '(y-x)$ and
$\{\phi_i(y),\phi_j(x)\}=-\delta_{ij}\delta '(y-x)$, then the
Poisson algebra induced on the $u_k's$ via the embedding
$u_k=u_k(\varphi_i,\phi_j)$ is $\W_{KP}^{(n-m)}$.}

{}From the point of view of the Lax equations, the first field $u_1$ in
$L= \d^n+u_1\d^{n-1}+u_2\d^{n-2}+...$ is a kind of spectator that can be
set to zero. However, from the point of view of the symplectic
structure, this constraint is second class, and the reduction changes
the Poisson bracket algebra. For example, performing this reduction in
$\GD_n$
yields the algebra $\W_n$. Correspondingly, in the case of a KP type
of Lax operator, $\W_{KP}^{(n)}$ reduces to the
non-linear $\hat{\W}^{(n)}_\infty$. The interest in
these reduced algebras
resides in the fact that they contain a Virasoro subalgebra
spanned by the
field $u_2$,  and that the higher generators $u_3,u_4,...$ can be
redefined to
transform like tensors of $Diff(S^1)$ under Poisson brackets with $u_2$.

The reduced Lax operador $\tilde L \equiv L|_{u_1=0}$ is in general
not consistent.
That is,
the Adler map \(adler) yields for generic $X$, another $J_{\tilde L}(X)$
 where
the coefficient of $\d^{n-1}$ is different from zero, and hence
``sticks out"
of the constrained manifold. This fact can be avoided if $X$ is
required
to satisfy the equation $res \comm{X}{\tilde L}=0$ \[Dickey].
Translated to the language of
Poisson brackets,
the reduced manifold ${\cal M}_{\tilde L} = \{u_1=0, u_2,u_3,..\}$
admits the same Poisson structure \(geldidos)
with the proviso that now the gradients are no more given by
\(gradients)
since in particular $\delta{F}/\delta{u_1}$ is ill defined. Instead they
should be defined by
$$
\gradtL F = \d^{-n}{\delta F\ov \delta u_1} +\sum_{i\geq 2} \d^{-n+i-1}
{\delta F\ov \delta u_i}
$$
with $\delta{F}/\delta{u_1}$ solved in terms
of $\delta F/\delta u_i,~(i=2,3,...)$
such that
$$
res\comm{\gradtL F}{\tilde L} = 0.
$$
holds identically.

It is now easy to check  that the above property is preserved
under the operation of taking the inverse of the Lax operator.
Explicitly, using the relation
$$
\gradtA F = - \tilde L \gradtL F \tilde L
$$
it follows that
$$
res\comm{\gradtA F}{\tilde A} = 0
 ~~~\Leftrightarrow ~~~~
res\comm{\gradtL F}{\tilde L} = 0.
$$

In the language of $\W$ algebras, it means that
if the fields in $a_i~ (i\geq 2)$
in $A=\d^n + a_2 \d^{n-2} + ...$ span $\W_n$,
the fields $u_i$ in $L = A^{-1}=
\d^{-n}+ u_2\d^{-n-2} + ...$ span ${ \hat{\W}}_{\infty}^{-n}$.

Unfortunately, we have not yet been able to reach
 a similar result in the context of
theorem I, but we hope to report on this point in the near future.

Summarizing, we have found that the  second
Gelfand-Dickey brackets exhibit a particularly simple
behavior under the product and
inversion of \pdo's, and this property is
crucial for inducing associated Poisson structures in the so called
constrained KP-hierarchies. To point other directions of further study,
it would be interesting
 to know to what extent the
Adler map is {\sl  completely determined} by this set of properties.

\ack
We would like to thank J.M. Figueroa O'Farrill
for enlightening conversations.
After this  note was written, he informed us
about a preprint by Yi Cheng
where  a similar approach is discussed.

\appendix{The Classical Limit}

This appendix is concerned with the classical
(also called dispersionless or commutative)
 limit of KP.  This dynamical system can also
be endowed with a Poisson structure
that arises as a classical limit of the
 second Gelfand-Dickey brackets \[Class].
The formal procedure requires first the mapping of
any ``normal ordered"\pdo\
(all $\d$'s to the right of the $u_i(x)$'s) to its
associated pseudo-differential symbol.
$$
A = \sum_{i}a_i \d^{n-i} \to   \hat A =\sum_{i}a_i \xi^{n-i}.
$$
The Poisson bracket of two functionals $F(\hat L)$
and $G(\hat L)$ is given by the
same expression \(compois) with $J_L(X)$ replaced by
$$
J^{\rm cl}_{\hat L}(\hat X) = [[\hat L,\hat X]]_+{\hat L}
 - [[{\hat L},(\hat X{\hat L})_+]] \(geldiclas)
$$
where for any two $\hat P$ and $\hat Q$, $\hat P\hat Q = \hat Q\hat P$
and $[[\hat P,\hat Q]]$ stands for the ordinary 2-dimensional
Poisson bracket, {\ie}
$$
[[\hat P,\hat Q]] = \pder{\hat P}{x}\pder{\hat Q}{\xi}
 -\pder{\hat Q}{x}\pder{\hat P}{\xi} .
$$
The following theorem is proven by a straightforward calculation
 of the same
type as before.

{\caps Theorem III}: {\sl Both
theorem I and II hold in the classical limit, \ie
 with $J_{\hat L}^{\rm cl}$ in place of $J_L$.}

In the classical case even a stronger result holds:
the mapping $L\to L^p$ where
$p\in\nats$ is a Poisson isomorphism. Namely let $u_i$ be the
fields in $L$, and $w_j$ those of $L^p$, and compute the Poisson
brackets among the $u_i$'s with the help of \(geldiclas).
The brackets induced on the $w_j(u_i)$ by the above mapping
coincide with those computed with
$J_{\hat L^p}^{\rm cl}$, up to a global factor
$p^2$ \[Class].

It is worth noting that even for the higher dimensional generalization
of the
classical Gelfand-Dickey brackets that were constructed in \[hidi] the
two previous theorems hold. This is not at all obvious since,
for example,
this construction uses not one but two different splittings.
 Again this fact calls
for an understanding of the {\sl uniqueness of the bracket under
 the requirement
that it ``replicates" upon multiplication of Lax operators}.

\refsout
\bye